\documentclass[preprint,review,12pt]{article}
\usepackage{amssymb, amsmath}
\usepackage{color}

\begin{document}

\begin{center}
\bigskip 

\textbf{STABILITY THRESHOLD FOR MULTIADDITIVE AND SYMMETRIC MAPPINGS}

Dan M. D\u{A}IANU

{\small \bigskip }
\end{center}

\textbf{Abstract}${\small .}${\small \ Z. Gajda showed that the control
functions of the form }$\epsilon \left( \left\Vert x\right\Vert
^{r}+\left\Vert y\right\Vert ^{r}\right) ${\small \ do not provide stability
for additive transformations if and only if }$r=1.${\small \ In this note we
prove a similar result for }$n${\small -additive and symmetric
functions.\bigskip }

${\small Keywords:}${\small \ }${\small n}${\small -additive mapings,
Hyers--Ulam--Rassias stability, stability threshold.}\bigskip

${\small MSC}$ ${\small (2010):}${\small \ 39B82; 39B52.}

\bigskip

{\large 1 \ \ Introduction}

The control functions of Hyers-Rassias type usually have a threshold of
stability, i.e. there is a family of control functions of this type for
which stability is flawed. The first and most famous example was given by
Gajda \cite{gaj91} for the stability of Cauchy's equation. Aoki \cite{ao50},
Rassias \cite{ra78}, for $r<1,$\ and then Gajda \cite{gaj91}, for $r>1$,
showed that if $S$ is a normed vectorial space, $B$ is a Banach space, $%
\epsilon >0$, and \textit{\ }$f:S\rightarrow B$\ is a function such that%
\textit{\ }%
\begin{equation}
\left\Vert f\left( x+y\right) -f\left( x\right) -f\left( y\right)
\right\Vert \leq \epsilon \left( \left\Vert x\right\Vert ^{r}+\left\Vert
y\right\Vert ^{r}\right) ,\text{ }x,y\in S  \label{1}
\end{equation}%
(with the convention that $\left\Vert 0\right\Vert ^{r}=1$ if $r\leq 0$)
then there exists a unique additive function \textit{\ }$a:S\rightarrow B$
such that 
\begin{equation*}
\left\Vert f\left( x\right) -a\left( x\right) \right\Vert \leq \epsilon 
\frac{2}{\left\vert 2-2^{r}\right\vert }\left\Vert x\right\Vert ^{r},\text{ }%
x\in S.
\end{equation*}

If $r=1,$ the assertion no longer remains valid. Let $f_{G}:%
\mathbb{R}
\rightarrow 
\mathbb{R}
$,%
\begin{equation*}
f_{G}\left( x\right) :=\sum\limits_{k=0}^{\infty }2^{-k}\zeta \left(
2^{k}x\right) ,x\in 
\mathbb{R}
,
\end{equation*}%
where the mapping $\zeta :%
\mathbb{R}
\rightarrow 
\mathbb{R}
$ is defined by

\begin{equation*}
\zeta \left( x\right) :=\left\{ 
\begin{array}{c}
\frac{\epsilon }{6},\text{ if\ }x\in \lbrack 1,\infty ), \\ 
\frac{\epsilon }{6}x,\text{ if }x\in (-1,1), \\ 
-\frac{\epsilon }{6},\text{ if\ }x\in (-\infty ,1].%
\end{array}%
\right.
\end{equation*}

\textbf{Lemma 1.1.} \cite{gaj91} \textit{The function }$f_{G}$ \textit{%
verifies\textit{\ (}\ref{1}), i.e. }$\left\vert f_{G}\left( x+y\right)
-f_{G}\left( x\right) -f_{G}\left( y\right) \right\vert $ $\leq \epsilon
\left( \left\vert x\right\vert +\left\vert y\right\vert \right) ,$ $x,y\in 
\mathbb{R}
$ \textit{but if }$\delta >0$\textit{\ and }$m:%
\mathbb{R}
\rightarrow 
\mathbb{R}
$\textit{\ is an additive mapping, there exists }$x_{m}\in 
\mathbb{R}
$\textit{\ such that }$\left\vert f_{G}\left( x_{m}\right) -m\left(
x_{m}\right) \right\vert >\delta \left\vert x_{m}\right\vert .$

\bigskip

We agree to say that $r=1$ is a \textit{stability}\ \textit{threshold }for
(Hyers -Ulam-Rassias stability of) additive mappings.

The\ symmetric and $n$-additive (or multiadditive) functions are important
tools in the characterizations of Fr\'{e}chet polynomials (see, for
instance, \cite{alm} for a new proof of the famous result of Fr\'{e}chet 
\cite{fr09}).

In this paper we complete the above-mentioned result and we find the
threshold of stability\textit{\ }for $n$-additive and symmetric mappings in
the particular case of a class of control functions of Hyers-Ulam-Rassias
type.

\bigskip

{\large 2 \ \ Main results}

In the following lines we consider that $S$ is an abelian semigroup, $B$ is
a Banach space, and $n$ is a positive integer.

For all symmetric function $g:S^{n}\rightarrow B,$ we denote by $%
D_{n}g:S^{n+1}\rightarrow B$ the mapping defined by $D_{n}g\left(
x_{1},x_{2}\right) :=$ $g\left( x_{1},x_{2}\right) -$ $g\left( x_{1}\right)
- $ $g\left( x_{2}\right) $ for $n=1,$ and for $n>1$ by 
\begin{eqnarray*}
D_{n}g\left( x_{1},...,x_{n+1}\right) &:&=\ g\left(
x_{1},...,x_{n-1},x_{n}+x_{n+1}\right) -g\left(
x_{1},...,x_{n-1},x_{n}\right) - \\
&&-g\left( x_{1},...,x_{n-1},x_{n+1}\right) .
\end{eqnarray*}%
We remember that the symmetric function $g:S^{n}\rightarrow B$ is $n$%
-additive if and only if $D_{n}g=0.$

If $\varphi :S^{n+1}\rightarrow \lbrack 0,\infty )$ is a function, $%
r_{n}\varphi :S^{n}\rightarrow \lbrack 0,\infty )$ denotes the mapping
defined by $r_{n}\varphi \left( x_{1}\right) :=\varphi \left(
x_{1},x_{1}\right) $ if $n=1,$ and for $n>1$%
\begin{equation*}
r_{n}\varphi \left( x_{1},...,x_{n}\right) :=\varphi \left(
2x_{1},...,2x_{n-1},x_{n},x_{n}\right) +2\varphi \left(
2x_{1},...,2x_{n-2},x_{n},x_{n-1},x_{n-1}\right) +
\end{equation*}%
\begin{equation*}
+\cdots +2^{n-2}\varphi \left(
2x_{1},x_{n},x_{n-1},...,x_{3},x_{2},x_{2}\right) +2^{n-1}\varphi \left(
x_{n},x_{n-1},...,x_{2},x_{1},x_{1}\right) .
\end{equation*}

Let $\varphi :S^{n+1}\rightarrow \lbrack 0,\infty ).$\ In \cite{dai} we
proved that if 
\begin{equation}
\sum\limits_{k=0}^{\infty }2^{-n\left( k+1\right) }\varphi _{n+1}\left(
2^{k}z\right) <\infty ,\text{ }z\in S^{n+1},  \label{2}
\end{equation}%
then 
\begin{equation*}
\widetilde{\varphi }^{+}:=\{\Phi :S^{n}\rightarrow \lbrack 0,\infty
)\left\vert {}\right. \Phi \left( y\right) \geq R_{n}^{+}\varphi \left(
y\right) \text{, and }\lim\limits_{k\rightarrow \infty }2^{-nk}\Phi \left(
2^{k}y\right) =0,\text{ }y\in S^{n}\}
\end{equation*}%
is a nonempty class, where 
\begin{equation*}
R_{n}^{+}\varphi \left( y\right) :=\sum\limits_{k=0}^{\infty }2^{-n\left(
k+1\right) }r_{n}\varphi \left( 2^{k}y\right) .
\end{equation*}

Also, in \cite{dai'} we proved that if $S$ is a commutative 2-divisible
commutative semigroup and%
\begin{equation}
\sum\limits_{k=0}^{\infty }2^{nk}\varphi _{n+1}\left( 2^{-k-1}z\right)
<\infty ,\text{ }z\in S^{n+1},  \label{3}
\end{equation}%
then 
\begin{equation*}
\widetilde{\varphi }^{-}:=\{\Phi :S^{n}\rightarrow \lbrack 0,\infty
)\left\vert {}\right. \Phi \left( y\right) \geq R_{n}^{-}\varphi \left(
y\right) \text{, and }\lim\limits_{k\rightarrow \infty }2^{nk}\Phi \left(
2^{-k}y\right) =0,\text{ }y\in S^{n}\}
\end{equation*}%
is a nontrivial class, where 
\begin{equation*}
R_{n}^{-}\varphi \left( y\right) :=\sum\limits_{k=0}^{\infty
}2^{nk}r_{n}\varphi \left( 2^{-k-1}y\right) .
\end{equation*}

Using the following elementary lemma, in \cite{dai} we have shown that the
functions which verifie (\ref{2}) constitute a class of control functions
that provide stability for $n$-additive and symmetric functions.

\bigskip

\textbf{Lemma 2.1. }\cite{dai} \textit{Let }$\left( b_{k}\right) _{k\in 
\mathbb{N}
}$\textit{\ be a sequence in }$B,$\textit{\ }$\left( \alpha _{k}\right)
_{k\in 
\mathbb{N}
}$\textit{\ be a sequence of positive numbers,\ and }$c>0$\textit{\ such
that }$\beta :=\sum\limits_{k=0}^{\infty }c^{-k-1}\alpha _{k}<\infty $%
\textit{\ and }$\left\Vert b_{k+1}-cb_{k}\right\Vert \leq \alpha _{k}$%
\textit{, for all }$k\in 
\mathbb{N}
$\textit{. Then }$\left( c^{-k}b_{k}\right) _{k\in 
\mathbb{N}
}$\textit{\ is a convergent sequence and }$\left\Vert b-b_{0}\right\Vert
\leq \beta $\textit{, where }$b=\lim\limits_{k\rightarrow \infty
}c^{-k}b_{k}.$

\bigskip

In the following lines we complete that result using control functions which
verifies (\ref{3}). For convenience, we reproduce from \cite{dai} the proof
of the result mentioned above also.

\bigskip

\textbf{Theorem 2.2.} \textit{Let }$\varphi :S^{n+1}\rightarrow \lbrack
0,\infty )$\textit{\ and }$g:S^{n}\rightarrow B$\textit{\ be a symmetric
function satisfying the inequality }%
\begin{equation}
\left\Vert D_{n}g\left( z\right) \right\Vert \leq \varphi \left( z\right) ,%
\text{ }z\in S^{n+1}.  \label{4}
\end{equation}

\textit{1. If \ }$\varphi $\textit{\ verifies \ (\ref{2}) then there exists
a unique symmetric }$n$\textit{-additive function }$a:S^{n}\rightarrow B$%
\textit{\ such that}%
\begin{equation}
\left\Vert \,g\left( y\right) -a\left( y\right) \right\Vert \leq \Phi \left(
y\right) ,\text{ }y\in S^{n}  \label{5}
\end{equation}%
\textit{for all }$\Phi \in \widetilde{\varphi }^{+}.$\textit{\ The symmetric 
}$n$\textit{-additive function }$a$\textit{\ is defined by }$a\left(
y\right) :=$\textit{\ }$\lim\limits_{k\rightarrow \infty }$\textit{\ }$%
2^{-nk}g\left( 2^{k}y\right) ,$\textit{\ }$y\in S^{n}.$

\textit{2. If }$S$\textit{\ is a }$2$\textit{-divisible abelian semigroup
and }$\varphi $\textit{\ verifies (\ref{3}), then there exists a unique
symmetric }$n$\textit{-additive function }$a:S^{n}\rightarrow B$\textit{\
which satisfies (\ref{5}) for all }$\Phi \in \widetilde{\varphi }^{-}.$%
\textit{\ The symmetric }$n$\textit{-additive function }$a$\textit{\ is
defined by }$a\left( y\right) :=$\textit{\ }$\lim\limits_{k\rightarrow
\infty }2^{nk}g\left( 2^{-k}y\right) ,$\textit{\ }$y\in \mathcal{S}^{n}.$

\textit{Proof.} Let $y=\left( x_{1},...,x_{n}\right) \in \mathcal{S}^{n}.$
Putting $z=\left( x_{1},...,x_{n-1},x_{n},x_{n}\right) $ in (\ref{5}), we get%
\newline
$\left\Vert D_{n}g\left( x_{1},...,x_{n-1},x_{n},x_{n}\right) \right\Vert
=\left\Vert g\left( x_{1},...,x_{n-1},2x_{n}\right) -\text{ }2g\left(
x_{1},...,x_{n}\right) \right\Vert \leq $\newline
$\hspace*{4.85cm}\leq \varphi \left( x_{1},...,x_{n-1},x_{n},x_{n}\right) $.

But $g$ is a symmetric function; therefore,

$\left\Vert g( 2y) \!-\!2^{n}g( y) \right\Vert \!=\!\left\Vert \right. [
g(2x_{1},2x_{2},...,2x_{n-1},2x_{n}) \!-\!2g(
2x_{1},2x_{2},...,2x_{n-1},x_{n})] \,+$

$2[g\left( 2x_{1},2x_{2},...,2x_{n-2},x_{n},2x_{n-1},\right) -2g\left(
2x_{1},2x_{2},...,2x_{n-2},x_{n},x_{n-1}\right) ]+\cdots $

$+2^{n-1}\left[ g\left( x_{n},x_{n-1},...,x_{2},2x_{1}\right) -\text{ }%
2g\left( x_{n},x_{n-1},...,x_{2},x_{1}\right) \right] \left. {}\right\Vert
\leq $

$\left\Vert D_{n}g\left( 2x_{1},...,2x_{n-1},x_{n},x_{n}\right) \right\Vert
+2\left\Vert D_{n}g\left( 2x_{1},...,2x_{n-2},x_{n},x_{n-1},x_{n-1}\right)
\right\Vert +\cdots $

$+2^{n-1}\left\Vert D_{n}g\left( x_{n},x_{n-1},...,x_{2},x_{1},x_{1}\right)
\right\Vert $ $\leq r_{n}\varphi \left( x_{1},...,x_{n}\right) ,$

hence,%
\begin{equation}
\left\Vert g\left( 2y\right) -2^{n}g\left( y\right) \right\Vert \leq
r_{n}\varphi \left( y\right) ,\text{ }y\in \mathcal{S}^{n}.  \label{6}
\end{equation}

1. Suppose that $\varphi $ verifies (\textit{\ref{2}) }and that $\Phi \in 
\widetilde{\varphi }^{+}$\textit{. }Replacing $y$ by $2^{k}y$ ($k\in 
\mathbb{N}
$) in (\textit{\ref{6}),} we get%
\begin{equation*}
\left\Vert \,g\left( 2^{k+1}y\right) -2^{n}g\left( 2^{k}y\right) \right\Vert
\leq r_{n}\varphi \left( 2^{k}y\right) ,\text{ }y\in \mathcal{S}^{n}.
\end{equation*}%
Applying Lemma 2.1 (for $b_{k}=g\left( 2^{k}y\right) ,$ $c=2^{n},\alpha
_{k}=r_{n}\varphi \left( 2^{k}y\right) $ and $\beta =R_{n}^{+}\varphi \left(
y\right) \leq \Phi \left( y\right) $), it follows that $\left(
2^{-nk}g(2^{k}y)\right) _{k\in 
\mathbb{N}
}$\ is a convergent\vspace*{0.03cm} sequence in $B,$\ and its limit, $a(y):=%
\mathit{\ }\lim\limits_{k\rightarrow \infty }\mathit{\ }2^{-nk}g\left(
2^{k}y\right) ,$ satisfies (\ref{5}). Since $g$ is symmetric, it follows
that $a$ is a symmetric function, too. From (\ref{4}) and (\ref{2})\ it
follows that%
\begin{equation*}
\lim\limits_{k\rightarrow \infty }\ 2^{-nk}\left\Vert
D_{n}g(2^{k}z)\right\Vert \leq \lim\limits_{k\rightarrow \infty }\
2^{-nk}\varphi \left( 2^{k}z\right) =0,
\end{equation*}%
whence $D_{n}a(z)=0,$ $z\in \mathcal{S}^{n+1}$, i.e. $a$ is a symmetric and $%
n$-additive mapping which satisfies (\ref{5}). If $a^{\prime }:\mathcal{S}%
^{n}\rightarrow B$ is an $n$-additive mapping and%
\begin{equation*}
\left\Vert g\left( y\right) -a^{\prime }\left( y\right) \right\Vert \leq
\Phi \left( y\right) ,\text{ }y\in \mathcal{S}^{n},
\end{equation*}%
since $a^{\prime }\left( 2^{k}y\right) =2^{nk}a^{\prime }\left( y\right) $,
and $\Phi \in \widetilde{\varphi }^{+},$\textit{\ }we have%
\begin{equation*}
\lim\limits_{k\rightarrow \infty }\left\Vert 2^{-nk}g\left( 2^{k}y\right)
-a^{\prime }\left( y\right) \right\Vert \leq \lim\limits_{k\rightarrow
\infty }2^{-nk}\in \Phi \left( 2^{k}y\right) =0,\text{ }y\in \mathcal{S}^{n},
\end{equation*}%
whence $a^{\prime }=a;$ therefore $a$ is the unique symmetric and $n$%
-additive mapping which satisfies (\ref{5}).

2. Suppose now that $\varphi $ verifies (\ref{3})\textit{\ }and that $\Phi
\in \widetilde{\varphi }^{-}$\textit{.} Replacing $y$ by $2^{-k-1}y$ ($k\in 
\mathbb{N}
$) in (\textit{\ref{6}\_,} we get%
\begin{equation*}
\left\Vert g\left( 2^{-k-1}y\right) -2^{-n}g\left( 2^{-k}y\right)
\right\Vert \leq 2^{-n}r_{n}\varphi \left( 2^{-k-1}y\right) ,\text{ }y\in 
\mathcal{S}^{n}.
\end{equation*}%
Using again Lemma 2.1 (for $b_{k}=g\left( 2^{-k}y\right) ,$ $c=2^{-n},$ $%
\alpha _{k}=2^{-n}r_{n}\varphi \left( 2^{-k-1}y\right) $ and $\beta
=R_{n}^{-}\varphi \left( y\right) \leq \Phi \left( y\right) $), it follows,
as in the first case, that $a\left( y\right) \!:=\!\lim\limits_{k\rightarrow
\infty }2^{nk}g\left( 2^{-k}y\right) $\vspace*{0.05cm} defines the unique
symmetric $n$-additive mapping which satisfies (\ref{5}). \ \ $\square $

The next consequence is a stability result in the Aoki-Rassias sense.

\bigskip

\textbf{Corollary 2.3. }Let $S$ be a normed space, $\epsilon >0$ and $r\neq
1.$ Suppose that $g:S^{n}\rightarrow B$ is a symmetric function such that%
\begin{equation}
\left\Vert D_{n}g\left( x_{1},...,x_{n+1}\right) \right\Vert \leq \epsilon
\left\Vert x_{1}\right\Vert ^{r}\cdots \left\Vert x_{n-1}\right\Vert
^{r}\left( \left\Vert x_{n}\right\Vert ^{r}+\left\Vert x_{n+1}\right\Vert
^{r}\right)  \label{7}
\end{equation}%
for all $x_{1},...,x_{n+1}\in S.$ Then there exists a unique $n$-additive
mapping $a:G^{n}\rightarrow B$\ for which%
\begin{equation}
\left\Vert g\left( x_{1},...,x_{n}\right) -a\left( x_{1},...,x_{n}\right)
\right\Vert \leq \epsilon \frac{2^{\left( n-1\right) \left( r-1\right) +1}}{%
\left\vert 2^{r}-2\right\vert }\left\Vert x_{1}\right\Vert ^{r}\cdots
\left\Vert x_{n}\right\Vert ^{r},  \label{8}
\end{equation}%
for all $x_{1},...,x_{n}\in S\backslash \left\{ 0\right\} $. If $r<1$\ then $%
a\left( y\right) :=$\ $\lim\limits_{k\rightarrow \infty }$\ $2^{-nk}g\left(
2^{k}y\right) ,$\ and if $r\in \left( 1,\infty \right) ,$ then $a\left(
y\right) :=$\ $\lim\limits_{k\rightarrow \infty }2^{nk}g\left(
2^{-k}y\right) ,$ $y\in S^{n}.$

\textit{Proof.} Let 
\begin{equation*}
\varphi \left( x_{1},...,x_{n+1}\right) :=\epsilon \left\Vert
x_{1}\right\Vert ^{r}\cdots \left\Vert x_{n-1}\right\Vert ^{r}\left(
\left\Vert x_{n}\right\Vert ^{r}+\left\Vert x_{n+1}\right\Vert ^{r}\right) .
\end{equation*}%
Then 
\begin{equation*}
r_{n}\varphi \left( x_{1},...,x_{n}\right) =\frac{2^{\left( n-1\right)
\left( r-1\right) +1}\left( 2^{nr}-2^{n}\right) }{2^{r}-2}.
\end{equation*}

1. Let $r<1.$ Then $\varphi $ verifies (\ref{2}). But, for $%
x_{1},...,x_{n}\neq 0,$ we have 
\begin{equation*}
R_{n}^{+}\varphi \left( x_{1},...,x_{n}\right) =\frac{r_{n}\varphi \left(
x_{1},...,x_{n}\right) }{2^{nr}-2^{n}}=\frac{2^{\left( n-1\right) \left(
r-1\right) +1}}{2^{r}-2}\left\Vert x_{1}\right\Vert ^{r}\cdots \left\Vert
x_{n}\right\Vert ^{r}.
\end{equation*}%
We apply Theorem 2.2 for $\Phi =R_{n}^{+}\varphi $ and we obtain (\ref{8})
for $a\left( y\right) =$\ $\lim\limits_{k\rightarrow \infty }$\ $%
2^{-nk}g\left( 2^{k}y\right) .$

2. Let $r>1.$ Then $\varphi $ verifies (\ref{3}) and 
\begin{equation*}
R_{n}^{-}\varphi \left( x_{1},...,x_{n}\right) =\frac{r_{n}\varphi \left(
x_{1},...,x_{n}\right) }{2^{n}-2^{nr}}=\frac{2^{\left( n-1\right) \left(
r-1\right) +1}}{2-2^{r}}\left\Vert x_{1}\right\Vert ^{r}\cdots \left\Vert
x_{n}\right\Vert ^{r}.
\end{equation*}%
We apply Theorem 2.2 for $\Phi =R_{n}^{-}\varphi $ and we obtain (\ref{8})
for $a\left( y\right) :=$\ $\lim\limits_{k\rightarrow \infty }2^{nk}g\left(
2^{-k}y\right) ,$ $y\in S^{n}.$ \ \ \ \ $\square $

\bigskip

\textbf{Stability threshold is }$r=1$

\bigskip

Let $S=B=%
\mathbb{R}
$ and $n\geq 2$. Let $\epsilon >0.$ The stability problem for $n$-additive
and symmetric mappings in the case $r=1$ is: there exists a positive
constant $\delta $ such that if $g:%
\mathbb{R}
^{n}\rightarrow 
\mathbb{R}
$ is a symmetric function and%
\begin{equation}
\left\vert D_{n}g\left( x_{1},...,x_{n},x_{n+1}\right) \right\vert \leq
\epsilon \left\vert x_{1}\right\vert \cdots \left\vert x_{n-1}\right\vert
\left( \left\vert x_{n}\right\vert +\left\vert x_{n+1}\right\vert \right) ,%
\text{ }x_{1},...,x_{n+1}\in 
\mathbb{R}
,\text{ }  \label{9}
\end{equation}%
there exists a unique symmetric and $n$-additive mapping $a:%
\mathbb{R}
^{n}\rightarrow 
\mathbb{R}
$ for which%
\begin{equation}
\left\vert g\left( x_{1},...,x_{n}\right) -a\left( x_{1},...,x_{n}\right)
\right\vert \leq \delta \left\vert x_{1}\right\vert \cdots \left\vert
x_{n}\right\vert ,\text{ }x_{1},...,x_{n}\in 
\mathbb{R}
.  \label{10}
\end{equation}

We give two examples. The first: a function $g$ which verifies (\ref{9}),
but for which there exist an infinity of symmetric and $n$-additive mappings
satisfying (\ref{10}). The second: a function $g$ which verifies (\ref{9}),
but for which there does not exist a symmetric and $n$-additive mapping $a$
satisfying (\ref{10}).

1. The symmetric function $g$ defined by $g\left( x_{1},...,x_{n}\right) :=%
\frac{\epsilon }{2}\left\vert x_{1}\right\vert \cdots \left\vert
x_{n}\right\vert $ verifies (\ref{9}), and, for all $\alpha \in \left[
-\delta +\frac{\epsilon }{2},\delta +\frac{\epsilon }{2}\right] $ $\ a\left(
x_{1},...,x_{n}\right) :=\alpha x_{1}\cdots x_{n}$ defines a symmetric $n$%
-additive mapping which satisfies (\ref{10}).

2. \ From Lemma 1.1, it follows that the function $f_{G}:%
\mathbb{R}
\rightarrow 
\mathbb{R}
$ verifies%
\begin{equation}
\left\vert f_{G}\left( x+y\right) -f_{G}\left( x\right) -f_{G}\left(
y\right) \right\vert \leq \epsilon \left( \left\vert x\right\vert
+\left\vert y\right\vert \right) ,\text{ }x,y\in 
\mathbb{R}
\label{11}
\end{equation}%
and, for all additive mapping $m:%
\mathbb{R}
\rightarrow 
\mathbb{R}
,$ there exists $x_{m}\in 
\mathbb{R}
$ for which%
\begin{equation}
\left\vert f_{G}\left( x_{m}\right) -a\left( x_{m}\right) \right\vert
>\delta \left\vert x_{m}\right\vert .  \label{12}
\end{equation}%
Let and $g:%
\mathbb{R}
^{n}\rightarrow 
\mathbb{R}
$\ be the symmetric function defined by%
\begin{equation*}
g\left( x_{1},...,x_{n}\right) :=f_{G}\left( x_{1}\right) x_{2}\cdots
x_{n}+x_{1}f_{G}\left( x_{2}\right) x_{3}\cdots x_{n}+\cdots +x_{1}\cdots
x_{n-1}f_{G}\left( x_{n}\right) .
\end{equation*}%
From (\ref{11}) it follows that $g$\ satisfies (\ref{9}). Suppose that $a:%
\mathbb{R}
^{n}\rightarrow 
\mathbb{R}
$ is a symmetric and $n$-additive mapping which verifies (\ref{10}). Then,
for $x_{1}=x_{2}=\cdots =x_{n-1}=1,$ and $x_{n}=x,$ we have from (\ref{10}):%
\begin{equation}
\left\vert f_{G}\left( x\right) -\left[ a\left( 1,...,1,x\right) -\left(
n-1\right) f_{G}\left( 1\right) x\right] \right\vert \leq \delta \left\vert
x\right\vert ,\text{ for all }x\in 
\mathbb{R}
.  \label{13}
\end{equation}%
But $m\left( x\right) :=a\left( 1,...,1,x\right) -\left( n-1\right)
f_{G}\left( 1\right) x$ defines an additive mapping and therefore (\ref{13})
contradicts (\ref{12}). \ \ $\square $

\begin{center}
\bigskip

{\small Dan M. D\u{a}ianu}

{\small Department of Mathematics, }

{\small "Politehnica" University of Timi\c{s}oara, }

{\small Victoriei Square, No.2, }

{\small 300006 Timi\c{s}oara, Romania, }

{\small dan.daianu@upt.ro }
\end{center}

\end{document}